\documentclass[journal]{IEEEtran}
\usepackage{amsfonts}
\usepackage{amsthm}
\usepackage{amsmath}
\usepackage{graphicx}
\usepackage{hyperref}

\begin{document}
\title{Corrections to ``Unified Laguerre polynomial-series-based
 distribution of small-scale fading envelopes''}
\author{Yin~Sun${}^\star$,~\IEEEmembership{Student Member,~IEEE,}
        \'Arp\'ad~Baricz,
        and~Shidong~Zhou${}^\star$,~\IEEEmembership{Member,~IEEE.}
\thanks{Copyright (c) 2010 IEEE. Personal use of this material is permitted.
However, permission to use this material for any other purposes must
be obtained from the IEEE by sending a request to
pubs-permissions@ieee.org.}
\thanks{${}^\star$Corresponding authors.}

\thanks{The work of Y. Sun and S. Zhou was supported by the National
Basic Research Program of China under Grant 2007CB310608, National
Natural Science Foundation of China under Grant 60832008 and Lab
project from Tsinghua National Lab on Information Science and
Technology (sub-project): Key technique for new distributed wireless
communications system.}
\thanks{The work of \'A. Baricz was supported by the
J\'anos Bolyai Research Scholarship of the Hungarian Academy of
Sciences and by the Romanian National Authority for Scientific
Research CNCSIS-UEFISCSU, project number PN-II-RU-PD 388/2010.}

\thanks{Y. Sun and S. Zhou are with the State Key Laboratory on Microwave
and Digital Communications, Tsinghua National Laboratory for
Information Science and Technology and Department of Electronic
Engineering, Tsinghua University, Beijing, China. Address: Room
4-407, FIT Building, Tsinghua University, Beijing 100084, People's
Republic of China. e-mail: sunyin02@gmail.com,
zhousd@tsinghua.edu.cn. \'A. Baricz is with the Department of
Economics, Babe\c{s}-Bolyai University, RO-400591, Cluj-Napoca,
Romania. e-mail: bariczocsi@yahoo.com.}
}

\maketitle
\begin{abstract}
In this correspondence, we point out two typographical errors in
Chai and Tjhung's paper and we offer the correct formula of the
unified Laguerre polynomial-series-based cumulative distribution
function (cdf) for small-scale fading distributions. A Laguerre
polynomial-series-based cdf formula for non-central chi-square
distribution is also provided as a special case of our unified cdf
result.
\end{abstract}
\IEEEpeerreviewmaketitle

\section{Two typographical errors in Chai
and Tjhung's paper \cite{Chai}}

In an interesting paper, Chai and Tjhung \cite{Chai} proposed new
unified probability density function (pdf) and cumulative
distribution function (cdf) formulas based on the generalized
Laguerre polynomial series expansion that cover a wide range of
small-scale fading distributions in wireless communications. Many
known Laguerre polynomial-series-based small-scale fading
distributions are special cases of their results, which include the
multiple-waves-plus-diffuse-power (MWDP) fading, $\kappa-\mu$
(noncentral chi), Nakagami$-m$, Rician (Nakagami$-n$), Nakagami$-q$
(Hoyt), Rayleigh, Weibull, and $\alpha-\mu$ (Stacy) distributions.

First we present two typographical errors of \cite{Chai}:

1. Chai and Tjhung stated that the orthogonality relation of
generalized Laguerre polynomial is expressed as \cite[Eq. (2)]{Chai}
\begin{eqnarray}
\int\limits_0^\infty e^{-x}L_l^\beta(x)L_n^\beta(x)dx =
\left\{\begin{array}{l}\Gamma(1+\beta){n+\beta\choose n},~~~\!
\textrm{if}~n=l,\\0,~~~~~~~~~~~~~~~~~~~\textrm{if}~n\neq
l,\end{array}\right.\!\!\!\!\!\!\nonumber
\end{eqnarray}
where $L_n^{\mu}$ is the generalized Laguerre polynomial, defined as
\cite[Eq. (22.3.9)]{Abramowitz}
\begin{eqnarray}\label{eq7}
\!\!\!\!&&L_n^{\mu}(x)=\frac{e^xx^{-\mu}}{n!}\frac{d^n}{dx^n}\left(e^{-x}x^{n+\mu}\right)=\sum_{k=0}^n{n+\mu\choose n-k}\frac{(-x)^k}{k!},\nonumber\\
\!\!\!\!&&~~~~~~~~~~~~~~~~~~~~~~~~~~~~~\mu>-1, n =
0,1,2\cdots\nonumber
\end{eqnarray}
However, there is a typographical error in \cite[Eq. (2)]{Chai}. The
correct formula is
\begin{eqnarray}
\int\limits_0^\infty x^\beta e^{-x}L_l^\beta(x)L_n^\beta(x)dx =
\left\{\!\!\begin{array}{l}\Gamma(1+\beta){n+\beta\choose n},~~~\!
\textrm{if}~n=l,\\0,~~~~~~~~~~~~~~~~~~~\textrm{if}~n\neq
l.\end{array}\right.\!\!\!\!\!\!\nonumber
\end{eqnarray}

2. The Laguerre polynomial expansion coefficients $C_n$ for the pdf
$f_X(x)$ of random variable $X$ was determined by \cite[Eq.
(10)]{Chai}. The fading envelop random variable $R$ is assumed to
satisfy
\begin{eqnarray}
x=r^\alpha/b,\nonumber
\end{eqnarray}
where $\alpha,b$ are parameters. Chai and Tjhung claimed that
\cite[Line -9 in p. 3990]{Chai}
\begin{eqnarray}
E_X\left[x^k\right]=E_R\left[r^{\alpha k}\right]/b,\nonumber
\end{eqnarray}
which contains a small typographical error. The correct expression
is given by
\begin{eqnarray} \label{eq3}
E_X\left[x^k\right]=E_R\left[r^{\alpha k}\right]/b^k.\nonumber
\end{eqnarray}
Using this and \cite[Eq. (10)]{Chai}, the coefficients $C_n$ can be
expressed by the $\alpha k$th moment of the fading envelope $R$.
%\begin{eqnarray}\label{eq3}
%C_n= \sum_{k=0}^n(-1)^k
%\frac{\Gamma(n+\beta+1)}{\Gamma(n-k+1)\Gamma(\beta+k+1)k!}\frac{E_R\left[r^{\alpha
%k}\right]}{b^k}.
%\end{eqnarray}
%The expression of $C_n$ can be simplified in some cases. The
%interested reader is referred to \cite[Section IV]{Chai}.

\section{the correct unified cdf formula}
Chai and Tjhung presented an unified cdf formula for the fading
envelope $R$ in \cite[Eq. (13)]{Chai}, i.e.,
\begin{eqnarray}\label{eq1}
F_R(R)\!\!\!\!\!\!\!\!\!\!&&=\int_0^Rf_R(r)dr \nonumber\\
&&=\frac{\alpha}{b^{\beta+1}}\sum_{n=0}^\infty
C_n\frac{n!}{\Gamma(n+\beta+1)}\nonumber\\
&&~~~~~~~~~\times\int_0^R\exp\left(-\frac{r^\alpha}{b}\right)
r^{\alpha(\beta+1)-1}L_n^\beta\left(\frac{r^\alpha}{b}\right)dr\nonumber\\
&&=\frac{1}{n}\left(\frac{R^\alpha}{b}\right)^{\beta+1}\exp\left(-\frac{R^\alpha}{b}\right)\nonumber\\
&&~~~~~~~~~\times \sum_{n=0}^\infty
C_n\frac{n!}{\Gamma(n+\beta+1)}L_{n-1}^{\beta+1}\left(\frac{R^\alpha}{b}\right),
\end{eqnarray}
where $\alpha,\beta,b$ are some parameters depending on the type of
fading distribution and are presented in \cite[Section V]{Chai} for
some special cases.

Unfortunately, the third line of the formula (\ref{eq1}) is not
valid, because $\frac{1}{n}$ is outside of the summation and the
generalized Laguerre polynomial $L_{n-1}^{\beta+1}$ is not defined
for $n=0$. The correct formula of the unified cdf is
\begin{eqnarray}\label{eq6}
F_R(R)\!\!\!\!\!\!\!\!\!\!&&=\left(\frac{R^\alpha}{b}\right)^{\beta+1}\exp\left(-\frac{R^\alpha}{b}\right)
\nonumber\\&&~~~~~~~~~~~\times\sum_{n=1}^\infty
C_n\frac{\Gamma(n)}{\Gamma(n+\beta+1)}L_{n-1}^{\beta+1}\left(\frac{R^\alpha}{b}\right)\nonumber\\
&&~~~~~~~~~~~+\frac{1}{\Gamma(\beta+1)}\gamma\left(\beta+1,\frac{R^\alpha}{b}\right),
\end{eqnarray}
where we used that $L_{0}^{\mu}(x)=1$ and $C_0=1.$ Here
$\gamma(\cdot,\cdot)$ is the lower incomplete gamma function,
defined by \cite[Eq. (6.5.2)]{Abramowitz}
\begin{eqnarray}\gamma(a,x)=\int_0^xt^{a-1}e^{-t}dt.\nonumber\end{eqnarray}
Therefore, the unified Laguerre polynomial-series-based cdf contains
one term expressed with the aid of the lower incomplete gamma
function.

\section{One special case of the unified cdf formula Eq. (\ref{eq6})}
We now provide a Laguerre polynomial-series-based cdf formula for
non-central chi-square distribution as a special case of Eq.
(\ref{eq6}). The pdf of non-central chi-square distribution
$\chi^2_{\nu,\lambda}$ with $\nu$ degrees of freedom and
non-centrality parameter $\lambda$ is given by \cite[Eq.
(29.4)]{Johnson}
\begin{eqnarray}\label{eq2}
f_{\chi^2_{\nu,\lambda}}(r)=
\frac{e^{-(r+\lambda)/2}}{2}\left(\frac{r}{\lambda}\right)^{\nu/4-1/2}I_{\nu/2-1}\left(\sqrt{\lambda
r}\right),\nonumber
\end{eqnarray} which is
equivalent to the pdf of $\kappa-\mu$ power distribution defined by
\cite[Eq. (2)]{Yacoub}. Tiku proposed the Laguerre polynomial
expansion of $f_{\chi^2_{\nu,\lambda}}(r)$ in 1965 \cite{Tiku},
which is given by (see also \cite[Eq. (29.11)]{Johnson})
\begin{eqnarray}\label{eq4}
&&\!\!\!\!\!\!\!\!\!\!\!\!\!\!\!\!\!\!f_{\chi^2_{\nu,\lambda}}(r) =
\frac{e^{-r/2}}{2}\left(\frac{1}{2}r\right)^{\left(\frac{\nu}{2}-1\right)}\sum_{j=0}^\infty
\frac{\left(-\frac{\lambda}{2}\right)^j}{\Gamma\left(\frac{1}{2}\nu+j\right)}L_j^
{\frac{\nu}{2}-1}\left(\frac{1}{2}r\right),\!\!\!\!\!\!\!\nonumber\\
&&~~~~~~~~~~~~~~~~~~~~~~~~~~~~~~~~~~~~~~~~~~\nu>0,\lambda\geq0.
\end{eqnarray}
Gideon and Gurland obtained a Laguerre polynomial expansion of the
cdf of non-central chi-square distribution in 1977 \cite{Gideon}
(see also \cite[Eq. (29.13)]{Johnson} with the parameter set of
$L_0$), i.e.,
\begin{eqnarray}\label{eq5}
F_{\chi^2_{\nu,\lambda}}(R) \!\!\!\!\!\!\!\!\!\!\!&&=
\int_0^Rf_{\chi^2_{\nu,\lambda}}(r)dr\nonumber\\
&& =
\left(\frac{R}{2}\right)^{\frac{\nu}{2}}\!\!\exp\left(-\frac{R}{2}\right)\sum_{n=1}^\infty
\frac{\left(-\frac{\lambda}{2}\right)^n}{\Gamma(n+\nu/2)n}
L_{n-1}^{\frac{\nu}{2}}\left(\frac{R}{2}\right)\!\!\nonumber\\
&&~~~~~~+\frac{1}{\Gamma(\nu/2)}\gamma\left(\frac{\nu}{2},\frac{R}{2}\right),~~~~\nu>0,\lambda\geq0.
\end{eqnarray}
Some other Laguerre polynomial expansions of
$F_{\chi^2_{\nu,\lambda}}(R)$ are available in \cite{Gideon} and
\cite[Eq. (29.13)]{Johnson}. We note that the pdf and cdf formulas
in Eq. (\ref{eq4}) and Eq. (\ref{eq5}) can be derived by choosing
$\alpha=1$, $b=2$, $\beta=\nu/2-1$ and $C_n=(-\lambda/2)^n/n!$ in
\cite[Eq. (11)]{Chai} and Eq. (\ref{eq6}) in previous section.

\end{document}